\documentclass[12pt]{article}
\usepackage{latexsym}

\def\bs{\begin{subequations}}
\def\es{\end{subequations}}

\catcode`\@=11

\newtoks\@stequation
\def\subequations{\refstepcounter{equation}
  \edef\@savedequation{\the\c@equation}%
  \@stequation=\expandafter{\theequation}
  \edef\@savedtheequation{\the\@stequation}
  \edef\oldtheequation{\theequation}%
  \setcounter{equation}{0}%
  \def\theequation{\oldtheequation\alph{equation}}}

\def\endsubequations{\setcounter{equation}{\@savedequation}%
  \@stequation=\expandafter{\@savedtheequation}%
  \edef\theequation{\the\@stequation}\global\@ignoretrue}

\makeatletter
        \renewcommand{\theequation}{\thesection.\arabic{equation}}%
        \@addtoreset{equation}{section}%
\makeatother

\renewcommand{\thefootnote}{\fnsymbol{footnote}}

\begin{document}

\begin{titlepage}
revised September 23, 2012 
 
submitted to International Journal of Modern Physics C

originally posted as arXiv:1209.1547

\begin{center}

{\large \bf Numerical Calculation of Bessel Functions \\}

\vskip 0.2in

Charles Schwartz\footnote{E-mail: schwartz@physics.berkeley.edu}

{\em Department of Physics,
     University of California\\
     Berkeley, California 94720}
        
\end{center}

\vskip .3in

\vfill

\begin{abstract}
A new computational procedure is offered to provide simple, accurate 
and flexible methods for using modern computers to give numerical 
evaluations of the various Bessel functions. The Trapezoidal Rule, 
applied to suitable integral representations, may become the method 
of choice for evaluation of the many Special Functions of 
mathematical physics.

\end{abstract}

\vfill

\end{titlepage}

\renewcommand{\thefootnote}{\arabic{footnote}}
\setcounter{footnote}{0}
\renewcommand{\thepage}{\arabic{page}}
\setcounter{page}{1}

\section{$K_{\nu} (z)$}

Let me start with this nice  integral 
representation for one type of Bessel function,
\begin{equation}
K_{\nu}(z) = \int_{0}^{\infty}dt\; cosh(\nu t) \;e^{-z\;cosh 
t}.\label{a1}
\end{equation}
See Appendix A for reference.

What we offer is a particular technique of numerical integration that 
is simple and rapidly convergent for infinite integrals of smooth 
functions: the Trapezoidal Rule, wisely applied.  The general theorem \cite{CS}
says the following.
\begin{eqnarray}
\sum_{n=-\infty}^{\infty}\;h\;f(nh) = \int_{-\infty}^{\infty} 
\;dx\;f(x) + {\cal{E}} \label{a2}\\
{\cal{E}} = \sum_{m=1}^{\infty}[F(2\pi m/h) + F(-2\pi m/h)],
\end{eqnarray}
where $F$ is the Fourier transform of the function $f$. If $f(x)$ is 
quite smooth, then its Fourier transform will decrease rapidly as the 
parameter h decreases. 

Another way of describing this method is that 
we see all of the correction terms in the Euler-Maclaurin 
series vanish. This means that as the interval $h$ decreases, the error 
in the Trapezoidal Rule decreases faster than any power of $h$. Precisely 
what that formula for the convergence rate is depends on the detailed 
analytic behavior of the function $f(x)$. Typically, one sees 
some kind of exponential decrease: ${\cal{E}} \sim e^{-\alpha 
h^{-\beta}}$. I do not offer any general theory about this formula, 
although particular models have been studied in the 
past.\cite{CS}\cite{JW}  
For example, if $\beta = 1$, one sees the number of correct decimal places 
double as 
one halves the interval $h$.  An essential virtue of this overall 
method is that the practitioner can see what the convergence looks 
like and can decide when to extend the calculation for better accuracy or to quit.
The examples shown in the Tables below offer such displays of rapidly 
converging output data.

For this method to be practical, one also needs the 
integrand to decrease rapidly as the integration variable goes toward 
infinity; and the sum over mesh points is just cut off when the 
contributions drop below the desired accuracy. This is sometimes 
helped by accelerating the rate of decay through a change of 
integration variable; and we shall show an example of this technique 
in Section 3. 

For readers who are 
new to this method, I recommend the simple exercise of computing the 
infinite Gaussian 
integral, $\int dx\;e^{-x^{2}}$, using the Trapezoidal rule (\ref{a2}). Here, the 
error formula turns out to be ${\cal{E}} \sim e^{-\pi^{2}/h^{2}}$. A recent 
mathematical review of this general technique may be found here \cite{JW}.

Applying this numerical technique directly to Eq. (\ref{a1}), we obtain the 
results shown in Tables 1a and 1b.  This looks quite good: rapid 
convergence to high accuracy at modest cost. The simple code for these 
calculations is shown in 
Appendix B; arbitrary (real, positive) values of z and $\nu$ may be input.

\vskip 0.5cm
Table 1a. Computations of Equation (\ref{a1}) using Method (\ref{a2})

\begin{tabular} {||c|c|c|c||}\hline
1/h & $K_{0}(z=0.1)$ & $K_{0}(z=1.0) \times 10$ & $K_{0}(z=10.) \times 10^{5}$ \\ \hline
1 &  2.42704 11398 56250 & 4.20936 51061 48591 & 2.28987 96730 52002\\ \hline
2 &  2.42706 90280 99576 & 4.21024 43651 11141 & 1.77845 16875 44865\\ \hline
4 &  2.42706 90247 02016 & 4.21024 43824 07083 & 1.77800 62316 16017\\ \hline
8 &  2.42706 90247 02016 & 4.21024 43824 07082 & 1.77800 62316 16764\\ \hline
16&  2.42706 90247 02016 & 4.21024 43824 07086 & 1.77800 62316 16764\\ \hline
\end{tabular}
\vskip 0.5cm
The number of mesh points used for the last line of data in Table 1a was 109, 73, 
39, respectively.

\vskip 0.5cm
Table 1b. Computations of Equation (\ref{a1}) using Method (\ref{a2})

\begin{tabular} {||c|c|c|c||}\hline
1/h & $K_{2.718}(z=0.01)\times 10^{-6}$ & $K_{2.718}(z=1.0) $ & 
$K_{2.718}(z=100.) \times 10^{45}$ \\ \hline
1 &  1.39714 10533 21390 & 4.54996 20838 02887 &  \\ \hline
2 &  1.40690 20983 29460 & 4.49904 64843 96175 & 9.30030 05343 36706\\ \hline
4 &  1.40690 07287 78440 & 4.49903 44319 18784 & 5.14859 62786 90992\\ \hline
8 &  1.40690 07287 78468 & 4.49903 44319 18744 & 4.83095 95172 64883\\ \hline
16&  1.40690 07287 78469 & 4.49903 44319 18749 & 4.83095 57412 19501\\ \hline
32 & & & 4.83095 57412 19519 \\ \hline
\end{tabular}
\vskip 0.5cm

The number of mesh points used for the last line of data in Table 1b 
was 150, 76, 
31, respectively.

\section{ $J_{\nu}(x)\; and \;N_{\nu}(x)$}
To get at the other Bessel functions we need to move our variables 
into the complex plane. Here is one standard relation \cite{MO}:
\begin{equation}
H_{\nu}^{(1)}(x) = J_{\nu}(x) + i N_{\nu}(x) =
\frac{2}{i \pi }\;e^{-i\nu \pi/2}\;K_{\nu}(-ix),\label{b1}
\end{equation}
and for now I'll take $x$ as a real variable. This function 
$N_{\nu}(x)$ is called $Y_{\nu}(x)$ by some authors.

Using the integral representation (\ref{a1})  now would not work well 
because what was previously a decaying exponential function of $t$ is 
now an oscillatory function. However, we can move the contour of 
integration for that variable $t$ so that it goes toward the line 
$Im(t) = \pi/2$; and this will restore the real exponential decay of 
the integrand at large distances.  Here is one prescription:
\begin{equation}
t = sinh(u) + \frac{i \pi}{2}\;tanh(u),\label{b2}
\end{equation}
and $u$ is a new real integration variable which we will treat 
according to the general method of Eq. (\ref{a2}). For an alternative 
approach, see Appendix C.

The simple program written for calculating $K_{\nu}(x)$ is now 
expanded to accommodate complex variables; this runs to about 50 lines 
of code in C and some numerical results 
are shown in Tables 2 and 3. The convergence looks quite good.
\vskip 0.5cm
Table 2. Calculation of $J_{\nu}(x)$ from Equation (\ref{b1}) using (\ref{b2}) in 
(\ref{a1}).

\begin{tabular} {||c|c|c|c||} \hline
1/h & $J_{1}(x=0.1) \times 10^{2}$ & $J_{1}(x=1.0) \times 10$  & $J_{1}(x=10.0)
\times 10^{2}$\\ \hline
1  &31.34519 12483 38983 & 4.84678 01345 03115 &-0.49254 74998 14282\\ \hline
2  &10.04145 83361 50352 & 4.40211 90106 01766 & 1.00684 96120 06995\\ \hline
4  & 4.66764 17956 85866 & 4.40051 65097 30195 & 4.31962 34829 07726\\ \hline
8  & 4.99456 06293 39347 & 4.40050 58776 70964 & 4.34727 46212 95072\\ \hline
16 & 4.99375 25888 30283 & 4.40050 58574 49333 & 4.34727 46168 86134\\ \hline
32 & 4.99375 26036 24231 & 4.40050 58574 49336 & 4.34727 46168 86136\\ \hline
64 & 4.99375 26036 24215 & 4.40050 58574 49336 & 4.34727 46168 86136\\ \hline 
\end{tabular}

\vskip 0.5cm

Table 3. Calculation of $N_{\nu}(x)$ from Equation (\ref{b1}) using (\ref{b2}) in 
(\ref{a1}).

\begin{tabular} {||c|c|c|c||} \hline
1/h & $N_{1}(x=0.1)$ & $N_{1}(x=1.0) \times 10$  & $N_{1}(x=10.0) 
\times 10$\\ \hline
1 & -8.81445 14807 36515 & -8.76505 96245 40165 & 5.92701 15605 77872\\ \hline
2 & -6.94259 74076 35009 & -7.79957 53906 29861 & 3.00580 17631 63178\\ \hline
4 & -6.44231 94399 89834 & -7.81226 14661 84539 & 2.48848 46077 69873\\ \hline
8 & -6.45895 10404 44470 & -7.81212 82132 14771 & 2.49015 42425 77341\\ \hline
16& -6.45895 10946 34644 & -7.81212 82130 02891 & 2.49015 42420 69539\\ \hline
32& -6.45895 10947 02030 & -7.81212 82130 02888 & 2.49015 42420 69539\\ \hline
64& -6.45895 10947 02026 & -7.81212 82130 02890 & 2.49015 42420 69538\\ \hline
\end{tabular}

\vskip 0.5cm

The number of mesh points used for the last line of data in these Tables was 
170, 143, 100, respectively.

\section{$I_{\nu}(x)$}

One can say that this Bessel function can almost always be 
effectively evaluated by using the power series. But, to complete the 
program offered in this paper, we should give an approach through an integral 
representation. Here is one, over an infinite range, which is found in the 
standard references. \cite{NIST}

\begin{equation}
I_{\nu}(x) =\frac{1}{2\pi i}\;\int_{\infty -i\pi}^{\infty + i\pi}\; 
dt \; e^{\nu t} \; e^{x\cosh(t)},\label{c1}
\end{equation}
which we can represent by the contour
\begin{equation}
t = cosh (u) + c + i\pi\;tanh(u),\label{c2}
\end{equation}
with the variable $u$ treated by Eq(\ref{a2}). The constant $c$ is 
arbitrary but if we choose $c = -1 + sinh^{-1}(\nu/x)$, then the 
point of stationary phase occurs at $u=0$.

This works numerically but not quite as nicely as our previous examples: there is greater 
loss of accuracy due to the oscillations of the integrand, especially 
at small values of x.

So, let's look at an alternative integral representation for this 
Bessel function:
\begin{equation}
I_{\nu}(x) = \frac{(x/2)^{\nu}}{\sqrt{\pi}\;\Gamma(\nu + 1/2)}\; 
\int_{0}^{\pi}d\theta\;sin^{2\nu}\theta\; cosh (x cos\theta).\label{c3}
\end{equation}
This is a purely real integral, with no oscillations, and the factor out in front takes care 
of the behavior at small x. But, Is it amenable to the Trapezoidal 
rule, Eq. (\ref{a2}),  for effective numerical integration?  If the 
index $\nu$ is an integer, the answer is yes, because the integrand 
is a periodic function over the interval of integration.

For more general application, however, we change variables to map the 
integration onto the entire real line:
\begin{equation}
cos \theta = t = tanh (u),\label{c4}
\end{equation}
and then treat the integration over $u$ by the Trapezoidal rule. Some 
results of this calculation are shown in Table 4.

\vskip 0.5cm
Table 4. Calculation of Equation (\ref{c3}) using (\ref{c4}) and (\ref{a2})

\begin{tabular} {||c|c|c|c||}\hline
1/h & $I_{2}(z=0.01)\times 10^{5}$ & $I_{2}(z=1.0) \times 10$ & 
$I_{2}(z=100.) \times 10^{-42}$ \\ \hline
1 & 1.30643 43443 38658 & 1.38166 26095 67052 & 1.20677 58487 72873\\ \hline
2 & 1.25004 98297 92758 & 1.35743 37442 58898 & 1.05923 08945 34915\\ \hline
4 & 1.25001 04167 00829 & 1.35747 66976 67279 & 1.05238 50353 94902\\ \hline
8 & 1.25001 04166 99218 & 1.35747 66976 70383 & 1.05238 43193 24316\\ \hline
16& 1.25001 04166 99218 & 1.35747 66976 70383 & 1.05238 43193 24312\\ \hline
\end{tabular}
\vskip 0.5cm
The number of mesh points used for the last line of data in Table 4 was
153,154, 184, respectively. 

\vskip 0.5cm

In all cases one may try to accelerate the rate of 
decay of the integral by such further transformations as $u = 
sinh(v)$ or $u = v^{3}$. For the same accuracy as shown in Table 4, 
this technique reduces the number of mesh points needed by a factor 
of 2 or 3. For the data shown in previous Tables, this did not 
produce improvement.

All numerical results given in this paper were obtained by 
calculations in ordinary double precision (16 decimals).

For the general use of Eq. (\ref{c3}) one needs accurate evaluation of 
the Gamma function for arbitrary argument; and that particular 
topic is mentioned in the following section.

\section{Discussion}

There 
are standard libraries available for the numerical computation of Bessel functions, 
which rely on a variety of old techniques.\footnote{I have been able 
to look into the open source libraries GSL and Ceres; how these things are 
done by Mathematica is unknown to me.} They use power series for 
small argument, asymptotic series for large arguments, various 
schemes involving recurrence relations and interpolation for intermediate arguments. As far as 
I can tell, the
technique presented here, which seems to cover all those bases in one 
grand sweep,   is new and has not been 
recognized before, although the general principles behind this technique have been 
known for some time.

This technique  appears to be uniquely valuable on several 
accounts: it gives high accuracy in rapid time; it is flexibly 
applicable to all sorts of Bessel functions; it is simple to program, 
relying on the standard computer routines for efficient evaluation of 
ordinary trigonometric functions (with real or complex arguments). It 
also invites the programmer (mathematician, physicist, engineer, or 
student) to 
be completely in charge of the analytical/numerical process, rather than relying 
on a big black box provided by some remote experts.

There is room for further exploration of various contours of the 
integration, since the simple ones used above may not be the most 
efficient. I expect that this general approach can be used for arbitrary complex 
values of the 
argument $z$ and order $\nu$, although I have not investigated those 
ideas. There are also some difficult problems in remote corners 
of Bessel function parameters \cite{JL}; and whether or not the present method 
might be better in all cases is an open question.

Finally, one can reasonably expect that other types of ``Special Functions'', 
beyond Bessel,  can be efficiently evaluated using these ideas if one has nice 
integral representations to start with.

For one example, the Gamma function can be nicely computed by using a 
contour integral:
\begin{equation}
1/\Gamma(z) = \frac{1}{2\pi i} \int\; dt \;e^{t}\;t^{-z}, \;\;\;\;\; 
t = 2-cosh(u) + i sinh(u),\label{d1}
\end{equation}
and the infinite integral over the real variable u is treated 
according to (\ref{a2}). See Table 5 for some numerical results.

\vskip 0.5cm
Table 5. Calculation of $\Gamma(z)$ using (\ref{d1}) and (\ref{a2}).

\begin{tabular}{||c|c|c|c||} \hline
1/h & $\Gamma(z=0.1)$ & $\Gamma(z=1+i10)$ \\ \hline
2 & 9.55495 13353 97527 &  \\ \hline
4 & 9.51350 56276 59931 &3.88422 34738 42538e-07+i1.13738 21528 
01491e-06\\ \hline
8 & 9.51350 76986 68703 &3.91892 92710 32289e-07+i1.12844 79696 
26640e-06\\ \hline
16 & 9.51350 76986 68744&3.91892 92708 81460e-07+i1.12844 79695 
84611e-06\\ \hline
32 & 9.51350 76986 68734&3.91892 92708 81405e-07+i1.12844 
79695 84617e-06 \\ \hline
64 & & 3.91892 92708 81394e-07+i1.12844 79695 84628e-06 \\ 
\hline
\end{tabular}
\vskip 0.5cm

Another famous function that is nicely handled this way is the zeta 
function.
\begin{equation}
\zeta(s) = \sum_{n=1}^{\infty}\; \frac{1}{n^{s}} = 
\frac{1}{2\Gamma(s)}\; \int_{-\infty}^{\infty}du\; e^{su}\; 
\frac{e^{-t/2}}{sinh(t/2)}, \;\;\;\; t=e^{u}.
\end{equation}
With a further acceleration by the change of variables $u = sinh(v)$, 
this gives efficient calculations with the Trapezoidal Rule.

\vskip 0.5cm
\setcounter{equation}{0}
\def\theequation{A.\arabic{equation}}
\boldmath
\noindent{\bf Appendix A}
\unboldmath
\vskip 0.5cm

Equation (\ref{a1}) may not be familiar. In fact, I do not find it in 
my standard reference book \cite{MO}. So, let me show that it is correct.  
Apply the 
differential operator:
\begin{eqnarray}
z^{2}[\frac{d^{2}}{dz^{2}} + \frac{1}{z}\frac{d}{dz} - 
\frac{\nu^{2}}{z^{2}} -1] K_{\nu}(z) = \\
\int_{0}^{\infty} dt \;cosh (\nu t) [z^{2}cosh^{2}t - z cosh t - 
\nu^{2}-z^{2}]e^{-z\;cosh t} = \\
\int_{0}^{\infty}\;dt \;cosh(\nu t) [\frac{d^{2}}{dt^{2}} - \nu^{2}] e^{-z\;cosh 
t},
\end{eqnarray}
and that equals zero. So this is some Bessel function. Now look at 
how this integral behaves as $z \rightarrow \infty$: It will be 
dominated by small values of $t$ and we readily find,
\begin{equation}
K_{\nu}(z) \rightarrow \sqrt{\pi/2z}\;e^{-z},
\end{equation}
which correctly identifies this particular Bessel function.

\vskip 0.5cm
\setcounter{equation}{0}
\def\theequation{B.\arabic{equation}}
\boldmath
\noindent{\bf Appendix B}
\unboldmath
\vskip 0.5cm
Here is the code used for calculaton of the data in Tables 1a and 1b.

\begin{verbatim}
#include <stdio.h>
#include <math.h>
int main (void)
{double z,t,term,sum,h,nu;
int n,k;
do {printf("Enter nu and z for Bessel function K_nu(z): ");
scanf("%lf %lf",&nu,&z);
h=2.0;
for(k=1;k<=8;k++) /* sequence of decreasing values of h */
{t=0.; /* the first point of the integral */
sum=0.5*exp(-z); 
n=1;
do {t+=h; /* successive points of the integral */
    term=cosh(nu*t)*exp(-z*cosh(t)); 
    sum+=term;
    n++;}
while(term/sum > 1.e-20); /* deciding when to quit */
printf("h= %lf, n=%d, ans = %24.16e\n",h, n, h*sum);
h=h/2.;}
} while(1);    
return 0;}
\end{verbatim}

\vskip 0.5cm
\setcounter{equation}{0}
\def\theequation{C.\arabic{equation}}
\boldmath
\noindent{\bf Appendix C}
\unboldmath
\vskip 0.5cm
Let me explore an alternative contour to the one given in Eq. 
(\ref{b2}). Move the contour integral of $t$ - which originally went 
out along the real axis from zero to infinity - as follows. First leg 
is $t = (0,0)$ to $t = (0,i\pi/2)$. Second leg is $t = (0, i\pi/2)$ to 
$t = (\infty, i\pi/2)$. 

It will simplify this discussion if I just take the special case of 
$\nu = 0$. The result of doing this leads to the following integral 
formulas.
\begin{eqnarray}
J_{0}(x) = \frac{2}{\pi}\int_{0}^{\pi/2}\;d \theta \; cos(x cos 
\theta),\label{B1} \\
N_{0}(x) =\frac{2}{\pi} \int_{0}^{\pi/2}\; d \theta \; sin(x cos 
\theta) - \frac{2}{\pi} 
\int_{0}^{\infty} \; ds \;e^{-x sinh(s)}.\label{B2}
\end{eqnarray}
These formulas are found in the standard literature \cite{MO}. But, 
are these integrals suitable for the particularly powerful method of 
numerical integration discussed in this paper? For $J_{0}$ the answer 
is yes; but this is only for special cases of the index $\nu$; and it 
works because the integrand is a periodic function. For 
$N_{0}$ the answer is no. Neither of the two integrals there 
satisfies the conditions for our method (although one could probably 
change the integration variables and make them conform to the type 
desired). An interesting question was whether the end-point 
correction terms from each of the two integrals at the pivot point $t 
= (0,i\pi/2)$ might cancel; and the answer appears to be negative. 
This is interpreted as emphasizing the role of continuous (analytic) 
functions and variables in getting the power of the present 
integration method. The contour discussed in this Appendix has a kink 
in it; and that is a spoiler.

\end{document}